\documentclass[epsf,12pt]{article}
\usepackage{a4}
\usepackage{latexsym}

\usepackage{epsfig}

\long\def\onefigure#1#2{
\begin{figure*}[tbp]
\begin{center}
#1
\end{center}
\caption{#2}
\end{figure*}
} 

\newcommand{\lipefig}[2]  
{\onefigure{\mbox{\psfig{file=#1.eps}}}{\label{f:#1} #2} }

\newtheorem{theorem}{Theorem}[section]

\newtheorem{prop}[theorem]{Proposition}

\newtheorem{lemma}[theorem]{Lemma}
\newtheorem{cl}[theorem]{Claim}
\newtheorem{cor}[theorem]{Corollary}

\newcommand{\owari}{\hfill\hbox{\rule{4pt}{8pt}}}
\def\qed{}
\long\outer\def\proof#1 #2\qed{\medbreak
  \noindent{{\bf Proof#1.}\hspace{\medskipamount}}{#2}\owari\par
\ifdim\lastskip<\medskipamount \removelastskip\penalty55\medskip\fi}
\def\pn{{\cal P}_n}
\def\area{{\rm Area}\,}
\def\per{{\rm Per}\,}
\def\inte{{\rm int}\,}
\def\conv{{\rm conv}\,}
\def\pos{{\rm pos}\,}
\def\dist{{\rm dist}\,}
\def\ck{{\cal K}}
\def\cs{{\cal S}}
\def\cf{{\cal F}}
\def\p{{\bf P}}
\def\r{{\bf R^2}}
\def\z{{\bf Z^2}}

\def\la{\lambda}

\def\de{\delta}
\def\al{\alpha}
\def\eps{\varepsilon}
\def\phi{\varphi}
\def\tri{\triangle}
\def\be{\beta}
\def\ga{\gamma}
\def\ga{\gamma}
\def\ga{\gamma}

\title{Jarn\'{\i}k's convex lattice $n$-gon for non-symmetric norms}
\author{ Imre B\'ar\'any, Nathana\"el Enriquez}

\begin{document}

\pagestyle {myheadings}

\markright{I. B\'ar\'any, N. Enriquez -  Jarn\'ik's convex lattice $n$-gon}

\maketitle
\begin{abstract}
What is the minimum perimeter of a convex lattice $n$-gon? This question was answered by
Jarn\'{\i}k in 1926. We solve the same question in the case when perimeter is measured by a
(not necessarily symmetric) norm.
\end{abstract}


\section{Introduction}
\setcounter{equation}{0} \setcounter{theorem}{0}

What is the minimal perimeter $L_n$ that a convex lattice polygon with $n$ vertices can have?
In 1926 Jarn\'{\i}k~\cite{jar} proved that $L_n=\frac
{\sqrt{6\pi}}9 n^{3/3}+O(n^{3/4})$. The aim of this paper is to extend this result to all,
not necessarily symmetric, norms in the plane. As usual, such a norm is defined by a convex
compact set $D \subset \r$ with $0 \in \inte D$, and the norm of $x \in \r$ is
\[
||x||=||x||_D= \min \{t\ge 0 : x \in tD\}.
\]

Let $\z$ be the lattice of integer points in $\r$, and write $\pn$ ($n\ge 3$) for
the set of all convex lattice $n$-gons in $\r$, that is, $P \in \pn$ if
$P=\conv\{z_1,\dots,z_n\}$ where $z_1,\dots,z_n \in \z$ are the vertices, in anticlockwise
order, of $P$. The $D$-perimeter of $P$ is defined by
\[
\per P= \per\!_D P=\sum_{i=1}^n ||z_{i+1}-z_i||_D
\]
where $z_{n+1}=z_1$ by convention. Note that for a non-symmetric $D$, $\per\!_D P$ depends on
the orientation of $P$ as well. Define now
\begin{equation}\label{eq:min}
L_n=L_n(D)= \min \{\per\!_D P: P \in \pn\}
\end{equation}
Since $D$ will be kept fixed throughout, we will often write $\per P$ and $L_n$ instead of
$\per_D P$ and $L_n(D)$.

In this paper we determine the asymptotic behaviour of $L_n(D)$ for all norms. We will also
show that, after suitable scaling, the minimizing polygons have a limiting shape. The same
results were proved by Maria Prodromou~\cite{pro} in 2005 in the case when $D$ is symmetric,
that is, $D=-D$. We will see that most of the difficulties in the non-symmetric case do not
come up in the symmetric one.

Define $\cf$ as the set of all positive continuous functions $r :[0,2\pi] \to {\bf R^+}$ with
$r(0)=r(2\pi)$. Such a function is the radial function of a {\sl starshaped} set in $\r$;
such a set contains the origin in its interior and the half-line starting at the origin in
direction $u(t)=(\cos t, \sin t)$ intersects its boundary at a single point which is at
distance $r(t)$ from the origin. We write $\cs$ for the set of all starshaped sets in $\r$.
Every convex compact set $K \subset \r$ with $0 \in \inte K$ is, of course, starshaped. We
denote by $\cf^c$ the set of radial functions of all such convex compact sets.

Let $r_0 \in \cf^c$ be the radial function of $D$.  The problem of
determining $L_n(D)$ is closely related to the following variational problem, to be denoted
by $VP(r_0)$. We seek a radial function $r\in \cf$ that minimizes
\begin{eqnarray}\label{eq:var}
&&\int_0^{2\pi} r^3(t)/r_0(t)dt \nonumber\\
\mbox{ subject to }&&\int_0^{2\pi} r^3(t)\cos t dt =0,\,\, \int_0^{2\pi} r^3(t)\sin
t dt=0,  \nonumber \\
\mbox{ and }&&\frac 12 \int_0^{2\pi}r^2(t)dt =1.
\end{eqnarray}

Assume $r(t)$ is the radial function of a convex (or starshaped) compact
set $K \subset \r$. Then the first condition says that the centre of
gravity, $g(K)$, of $K$ is at the origin, and the second condition says that $\area K =1$. We
will explain later the meaning of the function to be minimized. Using the results concerning
$L_n$ we will prove the following.

\begin{theorem}\label{th:var} There is a unique solution $r \in \cf$ to the variational problem. It
is the radial function of a convex compact set in $\r$ defined as the only function of the form
$${1\over r}={a\over r_0}+b\cos t +c \sin t$$
with $a>0$, $b,c\in {\bf R}$, that satisfies the constraints of $VP(r_0)$.
\end{theorem}

Notice that all the positive functions of the form ${a\over r_0}+b\cos t +c \sin t$ are
radial functions of a  convex set. Indeed, the sign of the curvature is given, in the
differentiable case,  by the sign of $({1\over r})"+{1\over r}$ which happens to be equal to
$a(({1\over r_0})"+{1\over r_0})$, which is always positive because
$D$ is convex. This result can easily be extended to the non differentiable case.

We mention further that the solution to $VP(r_0)$ is unique in a larger class than $\cf$. This will be clear from the proof. 

\section{Results and notations}\label{sec2}
\setcounter{equation}{0} \setcounter{theorem}{0}

Assume that the vertices of a minimizer $P_n \in \pn$ are $z_1,\dots,z_n$ in anticlockwise
order (which is the orientation giving the minimal $D$-perimeter). Then $E_n=\{z_2-z_1,\dots,
z_n-z_{n-1},z_1-z_n\}$ is the edge set of $P_n$. Define $C_n=\conv E_n$. Note that the $E_n$
determines $P_n$ uniquely (up to translation). Even more generally, the following is true.

\begin{prop} Suppose $V \subset \r$ is a finite set of vectors whose sum is zero. Assume further
that $u,v \in V$, $u=\la v$ with $\la>0$ implies that $u=v$. Then there is a unique (up to
translation) convex polygon whose edge set is equal to $V$.
\end{prop}

{\bf Proof.} This is very simple. One has to order (cyclically) the vectors in $V$ by
increasing slope as $v_1,\dots,v_n,v_1$. Then the polygonal path through the points
$0,v_1,v_1+v_2,v_1+v_2+v_3,\dots,v_1+\dots +v_n=0$ in this order is a convex polygon with
edge set $V$. Uniqueness is clear.\hfill$\Box$

We call this construction the {\sl increasing slope construction}. Here come our main
results. We let $\ck$ denote the family of all convex compact sets in $\r$ with non-empty
interior. For $K,L \in \ck$, $\dist (K,L)$ denotes their Hausdorff distance.

\begin{theorem}\label{th:onC} There is a unique $C \in \ck$ such that $\lim \dist(
(\area C_n)^{-1/2}C_n,C)=0$. Moreover, $g(C)=0$ and $\lim n^{-3/2}L_n(D)$ exists and equals
\[
\alpha (D)=\frac {\pi}{\sqrt 6}\int_C ||x||dx.
\]
\end{theorem}

We will prove the uniqueness part of Theorem~\ref{th:var} by showing that the radial function
of $C$ is the unique solution to the variational problem $VP(r_0)$.

\begin{theorem}\label{th:onPn} There is a convex set $P \subset \r$ such that the following holds.
Let $P_n$ be an arbitrary sequence of minimizers, of $L_n(D)$, translated so that $\min\{x:
(x,y) \in P_n\}$ is reached at the origin. Then $\lim \dist(n^{-3/2}P_n,P)=0$.
\end{theorem}

We explain in Section~\ref{sec:PndistP} how and why $P$ is determined uniquely by $C$. Moreover, it is shown in section 11 that the round shape found for $P$ in Jarnik's case is obtained if and only if the unit ball $D$ is given by an ellipse having a focus point at the origin.

To avoid some trivial complications in the proofs we assume that $D$ is strictly convex. We
emphasize however that the above results are valid without this extra condition. We make
another simplifying assumption, namely, that
\begin{equation}\label{eq:area}
\area D=1
\end{equation}
This is just a convenient scaling of the unit ball which leaves the set of minimizers, and
the corresponding $E_n$, $C_n$ and consequently $C,P$ unchanged.

The strategy of proof of the key Theorem 2.2 is as follows. We put together the following ingredients :
\begin{itemize}
\item almost all primitive vectors of $C_n$ belong to $E_n$ (Section 7),

\item the normalized convex hulls $(\area C_n)^{-1/2} C_n$ are sandwiched between two fixed Euclidean
balls (Section 6), so that the Blaschke selection theorem applies (Section 9),

\item the radial functions of the only possible limiting points of the sequence $(\area
C_n)^{-1/2} C_n$ are solutions of $VP(r_0)$ (Section 5). Moreover, the variational problem
$VP(r_0)$ has a unique solution (Section 8).
\end{itemize}

\section{Auxiliary lemmas}\label{aux}
\setcounter{equation}{0} \setcounter{theorem}{0}

We write $\p$ for the set of primitive vectors in $\z$, i.e., $z=(x,y) \in \z$ ($z\ne 0$) is
in $\p$ if $x$ and $y$ are relatively prime. The following two claims are very simple.

\begin{cl}\label{cl:decreas}For all $n \ge 3$, $L_n < L_{n+1}$.
\end{cl}

{\bf Proof.} Let $P_{n+1}=\conv\{z_0,z_1,\dots,z_n\}$ be a minimizer for $L_{n+1}$ and set
$P_n^*=\conv\{z_1,\dots,z_n\}$. Then $L_n \le \per P_n^* <L_{n+1}$. \hfill $\Box$

\begin{cl}\label{cl:primitive} $E_n \subset \p$.
\end{cl}

{\bf Proof.} Assume $P_n$ is a minimizer and the edge $z_2-z_1 \notin \p$, say. Then the
segment $[z_1,z_2]$ contains an integer $z \in \z$ distinct from $z_1,z_2$. The convex
lattice $n$-gon $\conv\{z_1,z,z_3,\dots,z_n\}$ has shorter $D$-perimeter than $P_n$ because
the triangle $\conv\{z_1,z_2,z_3\}$ contains the triangle $\conv\{z_1,z,z_3\}$ so the latter
has shorter $D$-perimeter. \hfill $\Box$

The following lemma will be useful when proving that most points in $C_n \cap \p$ belong to
$E_n$.

\begin{lemma}\label{la:simple} Assume $a,b \in E_n$ and $a \ne \pm b$. Let $T$ be the
parallelogram with vertices $0,a,b,a+b$. If $x,y \in (T \cap \p)\setminus E_n$ and $x\ne y$,
then $x+y \notin T$.
\end{lemma}

{\bf Proof.} If $x+y \in T$ were the case, then set $E^*=E_n \cup \{x,y,z\}\setminus \{a,b\}$
where $z=a+b-x-y$. The increasing slope construction works now because $\sum_{z \in E^*}z=0$
and gives rise to convex lattice $(n+1)$-gon $P$ if there is no $u \in E_n$ with $u=\la z$
with $\la>0$. If there is such a $u$, we replace $u$ and $z$ by $u+z$ in $E^*$, and the
increasing slope construction gives a convex lattice $n$-gon $P$. We claim
that $P$ has shorter $D$-perimeter than $P_n$. This clearly finishes the proof.

To prove $\per P < \per P_n$ we have to show that $\|x\|+\|y\|+\|z\| < \|a\|+\|b\|$. Assume
that the anticlockwise angle from $a$ to $b$ is smaller than $\pi$. Then $x,y,z \in
\pos\{a,b\}$ where $\pos\{a,b\}$ is the cone hull of $a$ and $b$. Order the vectors
$a,b,x,y,z$ by anticlockwise increasing slope. The outcome is $a,x,z,y,b$ say. Then the
triangle $\tri=\conv\{0,a,a+b\}$ contains the quadrilateral $Q=\conv\{0,x,x+z,x+y+z\}$ so the
latter has shorter $D$-perimeter. Now $a+b=x+y+z$ and
\[
\per Q=\|x\|+\|y\|+\|z\|+\|x+y+z\| < \per \tri = \|a\|+\|b\|+\|a+b\|,
\]
and $\per P < \per P_n$ follows.\hfill $\Box$
\bigskip
\bigskip

\ifx\JPicScale\undefined\def\JPicScale{1}\fi
\unitlength \JPicScale mm
\begin{picture}(146.88,61.25)(0,0)
\linethickness{0.5mm}
\put(30,10){\line(1,0){60}}
\linethickness{0.5mm}
\multiput(90,10)(0.12,0.36){125}{\line(0,1){0.36}}
\linethickness{0.1mm}
\put(45,55){\line(1,0){60}}
\linethickness{0.1mm}
\multiput(30,10)(0.12,0.36){125}{\line(0,1){0.36}}
\linethickness{0.1mm}
\multiput(30,10)(0.2,0.12){375}{\line(1,0){0.2}}
\linethickness{0.5mm}
\multiput(30,10)(0.6,0.12){63}{\line(1,0){0.6}}
\linethickness{0.1mm}
\multiput(30,10)(0.12,0.24){125}{\line(0,1){0.24}}
\linethickness{0.5mm}
\multiput(90,25)(0.12,0.24){125}{\line(0,1){0.24}}
\linethickness{0.5mm}
\multiput(67.5,17.5)(0.36,0.12){63}{\line(1,0){0.36}}
\linethickness{0.1mm}
\multiput(67.5,17.5)(0.12,0.24){125}{\line(0,1){0.24}}
\linethickness{0.1mm}
\multiput(45,40)(0.6,0.12){63}{\line(1,0){0.6}}
\linethickness{0.1mm}
\multiput(82.5,47.5)(0.36,0.12){63}{\line(1,0){0.36}}
\put(28.12,7.03){\makebox(0,0)[cc]{$\scriptstyle 0$}}

\put(52.5,34.84){\makebox(0,0)[cc]{}}

\put(38.44,48.91){\makebox(0,0)[cc]{}}

\put(111.56,35.78){\makebox(0,0)[cc]{}}

\put(146.88,51.25){\makebox(0,0)[cc]{}}

\put(90,55){\makebox(0,0)[cc]{}}

\put(91.25,52.12){\makebox(0,0)[cc]{$\scriptstyle z$}}

\put(122.5,38.12){\makebox(0,0)[cc]{}}

\put(10.62,61.25){\makebox(0,0)[cc]{}}

\put(116.88,57.5){\makebox(0,0)[cc]{}}

\put(108.88,57.38){\makebox(0,0)[cc]{$\scriptstyle a+b$}}

\put(114.12,54.38){\makebox(0,0)[cc]{$\scriptstyle =x+y+z$}}

\put(90.88,22.12){\makebox(0,0)[cc]{$\scriptstyle x+z$}}

\put(93.75,6.88){\makebox(0,0)[cc]{$\scriptstyle a$}}

\put(53.75,56.25){\makebox(0,0)[cc]{}}

\put(43.12,57.5){\makebox(0,0)[cc]{$\scriptstyle b$}}

\put(80,49.38){\makebox(0,0)[cc]{$\scriptstyle x+y$}}

\put(68.12,15){\makebox(0,0)[cc]{$\scriptstyle x$}}

\put(43.75,42.5){\makebox(0,0)[cc]{$\scriptstyle y$}}

\put(6.88,53.12){\makebox(0,0)[cc]{}}

\end{picture}

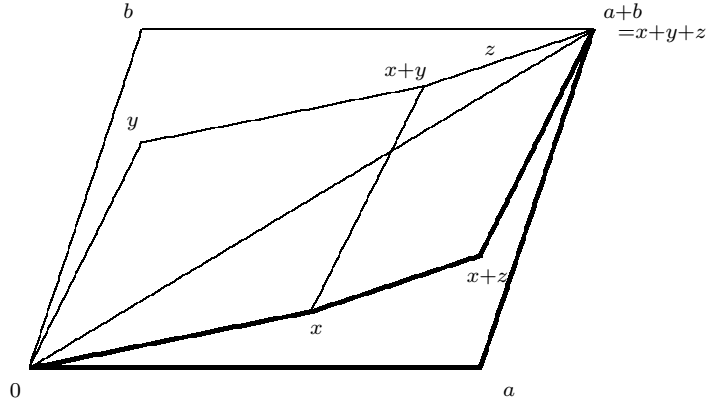
\centerline{Figure 1. The proof of Lemma 3.3}

\bigskip

We write $B$ for the Euclidean unit ball in $\r$ and $|x|$ for the Euclidean norm of $x \in
\r$. Since $D$ is compact convex and $0 \in \inte D$, there are positive constants $d_1,d_2$ such that
$d_1B \subset D \subset d_2B$, or, equivalently,
\[
d_1|x| \le \|x\| \le d_2|x|, \mbox{ for every } x \in D.
\]
In what follows $c,c_1,c_2,..$ denote positive constants independent of $n$. We will
also use Vinogradov's convenient $\ll$ notation: $f(n) \ll g(n)$  means that there are
positive constants $c$ and $n_0$ such that $cf(n) \le g(n)$ for all $n \ge n_0$. Of course,
the constants do not depend on $n$. But they depend on $D$, more precisely, they depend on
the constants $d_1,d_2$. $f(n) \gg
g(n)$ has the same meaning but with $f(n) \ge c g(n)$. We will also use the big Oh and little
oh notation.

We need some standard estimates on the distribution of lattice points and primitive points in
a convex body $K \in \ck$, see \cite{hawr} or \cite{bato} for a proof. Let $L$ denote the
Euclidean perimeter of $K$. We assume that $L >3$, say, but we think of $K$ as ``large''. In
fact, in most applications $L$ tends to infinity. The following estimate is simple and
well-known.
\begin{equation}\label{eq:z}
\left||K\cap \z| -\area K\right| \le  2L.
\end{equation}
This implies, with the standard method using the M\"obius function, that
\begin{equation}\label{eq:p}
\left||K\cap \p| -\frac 6{\pi^2}\area K\right| \le  3L \log L.
\end{equation}
Assume next that $f:\r \to {\bf R}$ is a $1$-homogeneous function, that is, $f(\la x)=\la
f(x)$ for every $x \in \r$ and $\la \ge 0$. Writing $M=\max\{|f(z)|: z\in K\}$ the following
estimates hold.
\begin{equation}\label{eq:fz}
\left| \sum _{z \in K \cap \z}f(z)- \int _K f(z)dz\right| \le 2 ML,
\end{equation}

\begin{equation}\label{eq:fp}
\left| \sum _{z \in K \cap \p}f(z)- \frac6{\pi^2}\int _K f(z)dz\right| \le 3ML \log L.
\end{equation}

The same estimates hold when $K$ is a (non-convex but) starshaped set whose boundary consists
of finitely many line segments. (Then, of course, the perimeter of $K$ is a finite number
$L$.) This fact will be needed in Section~\ref{sec:conn}.

These estimates will be used quite often in the case when $K=\la K_0$, 
and $\la \to \infty$ with $K_0$ fixed. Then formulae (\ref{eq:z}), 
(\ref{eq:p}), (\ref{eq:fz}), (\ref{eq:fp}) have the following simpler form:
\begin{equation}\label{eq:zla}
|K\cap \z| =\la^2\area K_0(1+O(\la^{-1})),
\end{equation}
\begin{equation}\label{eq:pla}
|K\cap \p| =\frac 6{\pi^2}\la^2 \area K_0((1+O(\la^{-1}\log \la)).
\end{equation}
\begin{equation}\label{eq:fzla}
\sum _{z \in K \cap \z}f(z)= \la^3\int _{K_0} f(z)dz+O(\la^{2}),
\end{equation}
\begin{equation}\label{eq:fpla}
\sum _{z \in K \cap \p}f(z)= \frac6{\pi^2}\la^3 \int _{K_0} f(z)dz+O(\la^{2} \log \la).
\end{equation}
The constant in the big Oh notation depends only on $K_0$. Here $K_0$ is either a convex set
or a starshaped set with boundary consisting of finitely many line segments.

\section{Bounding $L_n$}\label{sec:Ln}
\setcounter{equation}{0} \setcounter{theorem}{0}

In this section we give upper and lower bounds on $L_n$.

\begin{cl}\label{cl:Ln>} $L_n \gg n^{3/2}$.
\end{cl}

{\bf Proof.} Here we use the following {\sl density principle}. The sum of the lengths of $n$
distinct primitive vectors is at least as large as the sum of the lengths of the $n$ shortest
(distinct) primitive vectors. We will see the same principle in action a few more times.

Let $v_1,\dots,v_n$ be the $n$ shortest (in $D$-norm) vectors in $\p$ (ties broken
arbitrarily). Set $\la =\max\{\|v_i\|: i=1,\dots,n \}$. Then $(\inte \la D)\cap \p \subset
\{v_1,\dots,v_n\}\subset \la D$. The boundary of $\la D$ contains at most $\per\!_B \la D \le
2\pi d_2\la$ lattice points. So $|\la D \cap \p|-2\pi d_2 \la \le n \le |\la D \cap \p|$.
Using (\ref{eq:pla}) with $\la D$ (recalling $\area D=1$) gives
\[
|\la D\cap \p| =\frac 6{\pi^2} \la^2 (1+O(\la^{-1}\log \la)).
\]
This shows that $n=\frac {6}{\pi^2}\la^2(1 +O(\la^{-1}\log \la)$ implying that $\la=(\frac
{\pi}{\sqrt{6}}+o(1))n^{1/2}$. Using this in (\ref{eq:fpla}) with $\la D$ gives
\begin{eqnarray*}
L_n &\ge& \sum _1^n \|v_i\| \ge \sum _{z \in \inte (\la D) \cap \p}\|z \| \\
    &\ge& \left(\frac 6{\pi^2} -O(\la ^{-1}\log \la )\right)\la^3\int_D \|z\|dz \gg n^{3/2}.
\end{eqnarray*} \hfill $\Box$
\medskip

\begin{cl}\label{cl:Ln<} $L_n \ll n^{3/2}$.
\end{cl}

{\bf Proof.} Again, let $v_1,\dots,v_n$ be the $n$ shortest (in $D$-norm) vectors in $\p$ and
set $v_0=-\sum_1^n v_i$. By the increasing slope construction the vectors $v_0,v_1,\dots,v_n$
form the edge set of a unique (up to translation) convex lattice $n$-gon or $n+1$-gon. We
estimate its $D$-perimeter from above using the estimates on $\la$ from the previous proof.
\[
\sum _1^n \|v_i\| \le \sum _{z \in \la D}\|z \|
    \le \left(\frac 6{\pi^2} +O(\la^{-1} \log \la) \right)\la^3 \int_D
    \|z\|dz \ll n^{3/2},
\]
We need to estimate $\|v_0\|$ as well.
\[
\|v_0\| \ll |-v_0|= |v_0| \ll \|-v_0\| = \|\sum _1^n v_i\| \le \sum \|v_i\| \ll n^{3/2}.
\]
This shows that, indeed, $L_n \ll n^{3/2}$. \hfill $\Box$
\medskip

We mention that for a symmetric norm and for even $n$, the $n$ shortest vectors can be chosen in pairs $z,-z$ 
which is clearly optimal for $L_n$. The case of odd $n$ only causes only a minor difficulty.

\begin{cor}\label{cor:liminf} $\liminf n^{-3/2}L_n$ exists and equals $\al=\al(D)>0$, say.
\end{cor}

\section{Connection between $L_n$ and $VP(r_0)$}\label{sec:conn}
\setcounter{equation}{0} \setcounter{theorem}{0}

\begin{lemma}\label{l:conn} Assume $S \in \cs$ with $\area S=1$, $g(S)=0$. Then $r \in \cf$, the radial function of
$S$, is a feasible solution to $VP(r_0)$. Moreover, there is $Q_n \in \pn$ (for every $n\ge
3$) with
$$\lim n^{-3/2}\per Q_n=\frac {\pi}{\sqrt 6}\int_S\|z\|dz= \frac {\pi}{3\sqrt 6} \int_0^{2\pi} \frac {r^3(t)}{r_0(t)}dt.$$
\end{lemma}

We remark that the last identity follows from a simple integral transformation.
\smallskip

{\bf Proof.} Feasibility of $r$ is evident. We want to prove that for all $\eps>0$ (that we
will suppose small enough without restricting the generality),  there is $Q_n \in \pn$, for
every $n\ge 3$, with
$$\frac {\pi}{\sqrt 6}\int_S\|z\|dz-\eps\leq \liminf n^{-3/2}\per Q_n\leq\limsup n^{-3/2}\per Q_n\leq \frac {\pi}{\sqrt 6}\int_S\|z\|dz+\eps.$$
Since we would like to deal with sets whose perimeter can be defined and  controlled, we
introduce, for all $m\geq3$, the $m$-gon approximation of $S$, whose vertices are $r({2\pi
k\over m})u({2\pi k\over m})$ for $k=0,..., m-1$, recall that $u(t)=(\cos t,\sin t)$. The
sequence $S_m$ converges uniformly to $S$ as $m$ goes to infinity. Moreover, since $g(S)=0
\in \inte S$, there are constants $c_1,c_2>0$ such that, for $m$ large enough,
$$c_1B \subset S_m \subset c_2B.$$
We  fix now $m$ large enough so that the above condition is satisfied, as well as
$$\left\|{6\over\pi^2}\int_{S_m}zdz\right\|<c\eps$$
$$ \left|\frac {\pi}{\sqrt 6}\int_S\|z\|dz- \frac {\pi}{\sqrt 6}\int_{S_m}\|z\|dz\right|<c\eps$$
$$|\area S_m-1|<c\eps$$
where $c$ is a positive constant depending only on $S$ that will be adjusted later.

Now,  there is a minimal $\la>0$ (depending on $m$) so that $|\p \cap \la S_m|\ge n$. Let
$L_m$ denote the Euclidean perimeter of $S$. There are at most $\la L_m$ lattice points on
the boundary of $\la S_m$. Then, formula (\ref{eq:pla}) applies and shows  that
\[
|\p \cap \la S_m| =\left(\frac 6{\pi^2}\area S_m+O(\la^{-1}\log \la)\right) \la^2,
\]
implying $\la=\pi\sqrt{n/(6\area S_m)}(1+o(1))$.

 We apply formula (\ref{eq:fpla}) to $\la S$ with
$f(z)=z$, or more precisely with $f(z)=x$ and $f(z)=y$ where $z=(x,y)$ to get
\[
\sum_{z \in \p \cap \la S_m} z ={6\over\pi^2} \la^3 \int_{S_m}zdz+O(\la^2\log\la)
\]

Let $\p \cap \la S_m=\{z_1,\dots,z_l\}$ (of course $l \ge n$) and define $z_0=-\sum _1^l
z_i$. The previous equality implies that for $n$ large enough $\|z_0\| \leq 2c\eps\la^3$. The
increasing slope construction applies to $\{z_0,z_1,\dots,z_l\}$ and gives a convex lattice
$l$ or $l+1$-gon $T^n$. Note that $T^n$ has a {\sl special edge}, the one parallel to, and
having the same direction as, $z_0$. All other edges of $T^n$ are short, shorter than
$c_2d_2\la \ll n^{1/2}$ in $D$-norm. We claim now that, for a suitable choice of $c$
(depending only on $S$), and $\eps$ small enough,
$$\frac {\pi}{\sqrt 6}\int_S\|z\|dz-\eps\leq \liminf n^{-3/2}\per T_n\leq\limsup n^{-3/2}\per T^n \leq \frac {\pi}{\sqrt 6}\int_S\|z\|dz+\eps$$

We use again (\ref{eq:fpla}) this time with $f(z)=||z||$ to get
\[
\sum_1^l ||z_i|| = \frac 6{\pi^2}\la^3 \int_{\la S_m} ||z||dz (1+o(1))=\frac {\pi}{\sqrt 6}{n^{3/2}\over (\area S_m)^{3/2}} \int_{S_m} ||z||dz (1+o(1)).
\]
The claim follows since
$\per T^n$ differs from $\sum _1^m ||z_i||$ by $||z_0|| \leq 2c\eps\la^3$.

Finally, let $Q_n$ be the convex hull of $n$ consecutive vertices of $T^n$, including the two
endpoints of the special edge. Then $Q_n \in \pn$ and $\per Q_n \le \per T^n$ and also, $\per
Q_n$ is at least $\per T^n$ minus the sum of the $D$-length of the missing edges, which is
$\ll n$ as one can easily check. Thus $|\per Q_n -\per T^n| \ll n$. The requirements on the
constant $c$ are now clear. \hfill $\Box$
\medskip

We mention here that Lemma~\ref{l:conn} implies Claim~\ref{cl:Ln<} by simply choosing any $S
\in \cs$ with $g(S)=0$ and $\area S=1$, for instance the Euclidean disk centred at the origin and having area 1.

\section{Bounding $C_n$}\label{sec:Cn}

Our next target is to give bounds on the width and diameter of $C_n=\conv E_n$.

\begin{cl}\label{cl:wn>} The width of $E_n$, $w(E_n)$, satisfies $w(E_n) \gg n^{1/2}$.
\end{cl}

{\bf Proof.} Set $w=w(E_n)$. Clearly,
\[
L_n = \sum_{v \in E_n} \|v\| \gg \sum_{v \in E_n} |v|\ge  M_n(w),
\]
where $M_n(w)$ is the sum of the lengths of the $n$ shortest (in Euclidean norm) distinct
vectors in $\z$ lying in a strip of width $w$.

A simple yet technical computation, delayed to Appendix 1, shows that $w \le \ga n^{1/2}$
(where $\ga \in (0,1/2]$) implies $M_n(w) \gg n^{3/2}/\ga $. This finishes the proof of
Claim~\ref{cl:wn>}, because then $n^{3/2}\gg L_n \gg M_n(w) \gg n^{3/2}/\ga$ would lead to
contradiction if $\ga$ were too small.

\hfill$\Box$
\medskip

\begin{cl}\label{cl:diam<} Assume the smallest Euclidean ball centred at $0$ and
containing $E_n$ is $RB$. Then $R\ll n^{1/2}$.
\end{cl}

{\bf Proof.} Assume $a$ is the farthest point (in Euclidean distance) from the origin in
$E_n$. Then $|a|=R$. Claim~\ref{cl:Ln<} implies that $|a|\le L_n \ll n^{3/2}$. Since $w(E_n)
\gg n^{1/2}$ by the previous claim, there is a point $b \in E_n$ whose distance from the line
$\{x=ta: t\in {\bf R}\}$ is $\ge \frac 12 w(E_n) \gg n^{1/2}$.

The perimeter of the triangle $\tri=\conv\{0,a,b\}$ is $|a|+|b|+|a-b|\le 4|a|$ because
$|b|\le |a|$ and $|a-b|\le |a|+|b|\le 2|a|$. Here $\area \tri = \frac 12 |a|h$ where $h$ is
the corresponding height of $\tri$. Since $w(E_n)\ge n^{1/2}$, $h \gg n^{1/2}$.

Then by (\ref{eq:p}) for large enough $n$,
\[
\left| |\p \cap \frac 12 \tri| - \frac 6{\pi^2}\area \frac 12 \tri\right| \le 3\cdot 2|a|\log 2|a|
\ll h|a|\frac {\log |a|}{\sqrt n}\ll \area \tri \frac {\log n}{\sqrt n}
\]
implying that $|\p \cap \frac 12 \tri| \ge \frac 18 \area \tri$, again when $n$ is large enough.

Assume now that $\area \tri >16n$. Then $|\p \cap \frac 12 \tri| \ge 2n$. Since $|E_n| \le
n$, $\frac 12 \tri$ contains two distinct elements $x,y \in \p\backslash E_n$ and, evidently, $x+y \in
\tri$. Then $x,y,x+y \in \conv \{0,a,b,a+b\}$ contradicting Lemma~\ref{la:simple}.

Thus $\area \tri =\frac 12 |a|h\le 16n$, and so $R=|a| \ll n^{1/2}$. \hfill$\Box$
\bigskip

We need one more fact about $C_n$:

\begin{cl}\label{cl:inscrib} Assume $rB$ is the largest ball centred at $0$ and contained
in $C_n$. Then $r \gg n^{1/2}$.
\end{cl}

{\bf Proof.} Let $a$ be the nearest point to $0$ on the boundary of $C_n$. Thus $r=|a|$.
Define $E^+=E_n \cap \{x\in \r: ax>0\}$ and $E^-=E_n \cap \{x\in \r: ax<0\}$, and set
$f(x)=ax/|a|$ which is just the component of $x \in \r$ in direction $a$. To have simpler
notation we write $f(X)=\sum_{x\in X}f(x)$ when $X \in \r$ is a finite set. Since $\sum _{z\in
E_n}z=0$, $f(E^+) + f(E^-)=0$ (because $f(z)=0$ when $az=0$). We will show, however, that
$|a| \le \ga n^{1/2}$, for a suitably small $\ga>0$, implies that
\begin{equation}\label{eq:E+-}
f(E^+) + f(E^-)<0.
\end{equation}

Define $F^+=\{x\in \r: 0< f(x) \le \ga n^{1/2}\}\cap RB$ with $R \ll n^{1/2}$ from
Claim~\ref{cl:diam<}. The density principle tells now that $f(E^+)\le f(\p \cap F^+) \le f(\z \cap F^+)$ and the
last sum can be estimated as follows. Let $Q(z)$ be the unit cube centred at $z$. Again,
$\area Q(z) \cap F^+ \ge 1/4$ for all $z \in \z \cap F^+$. This implies that, for large
enough $n$,
\[
m:=|\z \cap F^+| \le \area F^+/4 \ll R|a| \ll \ga n.
\]
We use now (\ref{eq:fz}):
\[
\left|f(\z \cap F^+) -\int_{F^+}f(z)dz\right| \ll R|a|.
\]
It is easy to see that $\int_{F^+}f(z)dz \ll |a|^2R$ implying that $f(\z \cap F^+)\ll |a|^2R \ll \ga^2n^{3/2}$.

Define $F^-=\{x\in \r: 0> f(x) \ge -\la \ga n^{1/2}\}\cap RB$ where $\la >0$ is chosen so
that $F^-$ contains exactly $n-m-k$ lattice points. Here $k$ is the number of lattice points
on the line $ax=0$ so $k \le 2R+1 \ll n^{1/2}$. Note that $\la \ga n^{1/2}\ll R$ since $E_n \subset RB$ consists of exactly $n$ vectors.
Choosing $\ga$ small enough guarantees that
$m< 0.1n$ which, in turn, guarantees that $\la >1$ and further, that $|F^-\cap \z|\ge 0.8n$.
The Euclidean perimeter of $F^-$ is at most 
$4R+\la \ga n^{1/2}\ll R$ and (\ref{eq:z}) shows that $\left||F^-\cap \z|-\area F^-\right| \ll R$. 
Clearly $\area F^- \ll R\la\ga n^{1/2}$, implying that
\[
0.8n < |F^-\cap \z| \le \left(1+O\left( \frac {1}{\la\ga n^{1/2}}\right)\right)\area F^- \ll R \la \ga
n^{1/2}\ll \la \ga n,
\]
which implies $\la \ga \gg 1$.

The density principle says now that $f(E^-) \le f(F^-)$ (note that $f$ is negative on $F^-$ and $E^-$), and $f(F^0)$ can be estimated using (\ref{eq:fz}):
\[
\left|f(F^-)-\int_{F^-}f(z)dz\right| \ll R^2 \ll n,
\]
because $\max\{|f(x)|: x \in F^-\}\le R$. Now $f(z)$ is negative on $F^-$. It is easy to check
that $\la^2 \ga^2 nR \ll -\int_{F^-}f(z)dz \ll \la^2 \ga^2 nR$.
So we have
\[
-f(F^-)\ge \int_{F^-}-f(z)dz +O(n)\gg \int_{F^-}-f(z)dz \gg \la^2\ga^2 n R \gg n^{3/2}
\]
This shows that (\ref{eq:E+-}) indeed holds if $\ga>0$ is chosen small enough because
$0<f(\z\cap F^+) \ll \ga^2n^{3/2}$ and $-f(\z\cap F^-) \gg n^{3/2}$.
 \hfill$\Box$ \medskip

\begin{cor}\label{cor:sandw} There are positive numbers $r$ and $R$ (depending only on $D$) such that for all
$n\ge 3$
$$rB \subset (\area C_n)^{-1/2} C_n \subset RB.$$
\end{cor}

\section{Almost all primitive points of $C_n$ are in $E_n$}\label{CnandEn}

We begin by stating a geometric lemma which is about a special kind of approximation. The
technical proof is postponed to Appendix 2.

\begin{lemma}\label{le:approx} Assume $K \in \ck$ is a convex polygon with $rB \subset K \subset
RB$. Then for every $\de \in (0,0.02(r/R)^2]$ there are vertices $v_1,\dots,v_m$ of $K$ such
that with $Q=\conv\{v_1,\dots,v_m\}$ the following holds:
\begin{itemize}
\item $Q \subset K \subset (1+4R^2r^{-2}\de) Q$,
\item  for all $i$, the angle $\angle v_i0v_{i+1}$ is at least $\de$.
\end{itemize}
\end{lemma}

\begin{lemma}\label{le:CnandEn}For every $\eps >0$ there is $n_0=n_0(\eps,D)$ such that for all
$n \ge n_0$, $(1-\eps)C_n \cap \p \subset E_n$.
\end{lemma}

{\bf Proof.} Let $r_n$, resp. $R_n$ be the maximal, minimal radius such that $r_nB \subset
C_n \subset R_nB$. It follows from Claims~\ref{cl:wn>} and \ref{cl:inscrib} that $R_n/r_n \le
c$ with a suitable positive constant depending only on $D$. Thus Lemma~\ref{le:approx} can be
applied with $K=C_n$ and $\de= \eps/(8c^2)$ (if $\eps\le 0.02/8$ which we can clearly
assume). We get a polygon $Q=\conv\{v_1,\dots,v_m\}$ satisfying $C_n \subset (1+\eps/2)Q$.

Assume, contrary to the statement of the lemma, that there is an $x \in (1-\eps)C_n \cap \p
\setminus E_n$. One of the cones $\pos\{v_i,v_{i+1}\}$ contains $x$, say in the cone $W:=
\pos\{v_1,v_2\}$. Define $\tri=\conv\{0,v_1,v_2\}$. Thus $\tri \subset C_n \cap W \subset
(1+\eps/2)\tri$.
As $x \in (1-\eps)C_n\cap W$, $v_1+v_2-x \in W \setminus (1+\eps)\tri$. The
triangle $\tri^*= \left((v_1+v_2-x)-W\right)\setminus (1+\eps/2)\tri$ is disjoint from $C_n$.
We claim that it contains a primitive point $y$. This will finish the proof since then
$x,y,x+y$ all lie in the parallelogram with vertices $0,v_1,v_2,v_1+v_2$ contradicting
Lemma~\ref{la:simple}.

We prove the claim by using (\ref{eq:p}): $\area \tri^* \gg \eps^3 n$ because its angle at
$v_1+v_2-x$ is at least $\de$, and the neighbouring sides are of length at least
$\eps|v_1|/2$ and $\eps|v_2|/2$ and $|v_1|,|v_2|\gg n^{1/2}$. Further,
its perimeter is at most $|v_1|+|v_2|+|v_1-v_2| \ll n^{1/2}$. Thus
\[
\left||\tri^* \cap \p|-\frac 6{\pi^2}\area \tri^* \right| \ll (\log n) n^{1/2}.
\]
Here $\frac 6{\pi^2}\area \tri^*$ is of order $\eps^3 n$ and the error term is of
order $(\log n) n^{1/2}$. Since $\eps $ fixed, $\tri^*$ contains a primitive vector if $n$ is large
enough. \hfill $\Box$\bigskip

\bigskip
\ifx\JPicScale\undefined\def\JPicScale{1}\fi
\unitlength \JPicScale mm
\begin{picture}(120,54.84)(0,0)
\linethickness{0.1mm}
\put(19.53,9.84){\line(1,0){52.5}}
\multiput(19.53,9.84)(0.18,0.12){250}{\line(1,0){0.18}}
\multiput(72.03,9.84)(0.18,0.12){250}{\line(1,0){0.18}}
\put(64.53,39.84){\line(1,0){52.5}}
\linethickness{0.1mm}
\multiput(72.03,9.84)(0.16,0.47){1}{\line(0,1){0.47}}
\multiput(72.19,10.31)(0.15,0.47){1}{\line(0,1){0.47}}
\multiput(72.35,10.78)(0.15,0.47){1}{\line(0,1){0.47}}
\multiput(72.49,11.26)(0.14,0.48){1}{\line(0,1){0.48}}
\multiput(72.63,11.73)(0.13,0.48){1}{\line(0,1){0.48}}
\multiput(72.76,12.21)(0.12,0.48){1}{\line(0,1){0.48}}
\multiput(72.88,12.69)(0.11,0.48){1}{\line(0,1){0.48}}
\multiput(72.99,13.18)(0.1,0.49){1}{\line(0,1){0.49}}
\multiput(73.1,13.66)(0.09,0.49){1}{\line(0,1){0.49}}
\multiput(73.19,14.15)(0.09,0.49){1}{\line(0,1){0.49}}
\multiput(73.28,14.64)(0.08,0.49){1}{\line(0,1){0.49}}
\multiput(73.35,15.13)(0.07,0.49){1}{\line(0,1){0.49}}
\multiput(73.42,15.62)(0.06,0.49){1}{\line(0,1){0.49}}
\multiput(73.48,16.12)(0.05,0.49){1}{\line(0,1){0.49}}
\multiput(73.53,16.61)(0.04,0.49){1}{\line(0,1){0.49}}
\multiput(73.57,17.11)(0.03,0.5){1}{\line(0,1){0.5}}
\multiput(73.61,17.6)(0.02,0.5){1}{\line(0,1){0.5}}
\multiput(73.63,18.1)(0.02,0.5){1}{\line(0,1){0.5}}
\multiput(73.65,18.59)(0.01,0.5){1}{\line(0,1){0.5}}
\multiput(73.65,19.59)(0,-0.5){1}{\line(0,-1){0.5}}
\multiput(73.64,20.08)(0.01,-0.5){1}{\line(0,-1){0.5}}
\multiput(73.62,20.58)(0.02,-0.5){1}{\line(0,-1){0.5}}
\multiput(73.59,21.08)(0.03,-0.5){1}{\line(0,-1){0.5}}
\multiput(73.55,21.57)(0.04,-0.5){1}{\line(0,-1){0.5}}
\multiput(73.51,22.07)(0.05,-0.49){1}{\line(0,-1){0.49}}
\multiput(73.45,22.56)(0.06,-0.49){1}{\line(0,-1){0.49}}
\multiput(73.39,23.05)(0.06,-0.49){1}{\line(0,-1){0.49}}
\multiput(73.32,23.54)(0.07,-0.49){1}{\line(0,-1){0.49}}
\multiput(73.23,24.03)(0.08,-0.49){1}{\line(0,-1){0.49}}
\multiput(73.14,24.52)(0.09,-0.49){1}{\line(0,-1){0.49}}
\multiput(73.05,25.01)(0.1,-0.49){1}{\line(0,-1){0.49}}
\multiput(72.94,25.49)(0.11,-0.48){1}{\line(0,-1){0.48}}
\multiput(72.82,25.98)(0.12,-0.48){1}{\line(0,-1){0.48}}
\multiput(72.7,26.46)(0.12,-0.48){1}{\line(0,-1){0.48}}
\multiput(72.56,26.93)(0.13,-0.48){1}{\line(0,-1){0.48}}
\multiput(72.42,27.41)(0.14,-0.48){1}{\line(0,-1){0.48}}
\multiput(72.27,27.88)(0.15,-0.47){1}{\line(0,-1){0.47}}
\multiput(72.11,28.35)(0.16,-0.47){1}{\line(0,-1){0.47}}
\multiput(71.95,28.82)(0.17,-0.47){1}{\line(0,-1){0.47}}
\multiput(71.77,29.29)(0.18,-0.46){1}{\line(0,-1){0.46}}
\multiput(71.59,29.75)(0.09,-0.23){2}{\line(0,-1){0.23}}
\multiput(71.4,30.21)(0.1,-0.23){2}{\line(0,-1){0.23}}
\multiput(71.2,30.66)(0.1,-0.23){2}{\line(0,-1){0.23}}
\multiput(70.99,31.11)(0.1,-0.23){2}{\line(0,-1){0.23}}
\multiput(70.77,31.56)(0.11,-0.22){2}{\line(0,-1){0.22}}
\multiput(70.55,32)(0.11,-0.22){2}{\line(0,-1){0.22}}
\multiput(70.32,32.44)(0.12,-0.22){2}{\line(0,-1){0.22}}
\multiput(70.08,32.88)(0.12,-0.22){2}{\line(0,-1){0.22}}
\multiput(69.83,33.31)(0.12,-0.22){2}{\line(0,-1){0.22}}
\multiput(69.57,33.73)(0.13,-0.21){2}{\line(0,-1){0.21}}
\multiput(69.31,34.16)(0.13,-0.21){2}{\line(0,-1){0.21}}
\multiput(69.04,34.57)(0.13,-0.21){2}{\line(0,-1){0.21}}
\multiput(68.76,34.98)(0.14,-0.21){2}{\line(0,-1){0.21}}
\multiput(68.48,35.39)(0.14,-0.2){2}{\line(0,-1){0.2}}
\multiput(68.19,35.79)(0.15,-0.2){2}{\line(0,-1){0.2}}
\multiput(67.89,36.19)(0.15,-0.2){2}{\line(0,-1){0.2}}
\multiput(67.58,36.58)(0.1,-0.13){3}{\line(0,-1){0.13}}
\multiput(67.27,36.97)(0.1,-0.13){3}{\line(0,-1){0.13}}
\multiput(66.95,37.35)(0.11,-0.13){3}{\line(0,-1){0.13}}
\multiput(66.62,37.72)(0.11,-0.12){3}{\line(0,-1){0.12}}
\multiput(66.29,38.09)(0.11,-0.12){3}{\line(0,-1){0.12}}
\multiput(65.95,38.45)(0.11,-0.12){3}{\line(0,-1){0.12}}
\multiput(65.61,38.81)(0.12,-0.12){3}{\line(0,-1){0.12}}
\multiput(65.25,39.16)(0.12,-0.12){3}{\line(1,0){0.12}}
\multiput(64.89,39.5)(0.12,-0.11){3}{\line(1,0){0.12}}
\multiput(64.53,39.84)(0.12,-0.11){3}{\line(1,0){0.12}}
\multiput(64.53,39.84)(0.12,-0.48){62}{\line(0,-1){0.48}}
\linethickness{0.1mm}
\multiput(62.66,9.84)(0.1,0.23){2}{\line(0,1){0.23}}
\multiput(62.85,10.3)(0.09,0.23){2}{\line(0,1){0.23}}
\multiput(63.04,10.77)(0.17,0.47){1}{\line(0,1){0.47}}
\multiput(63.21,11.24)(0.16,0.47){1}{\line(0,1){0.47}}
\multiput(63.36,11.71)(0.15,0.48){1}{\line(0,1){0.48}}
\multiput(63.51,12.19)(0.14,0.48){1}{\line(0,1){0.48}}
\multiput(63.65,12.67)(0.12,0.48){1}{\line(0,1){0.48}}
\multiput(63.77,13.16)(0.11,0.49){1}{\line(0,1){0.49}}
\multiput(63.89,13.64)(0.1,0.49){1}{\line(0,1){0.49}}
\multiput(63.99,14.13)(0.09,0.49){1}{\line(0,1){0.49}}
\multiput(64.08,14.63)(0.08,0.49){1}{\line(0,1){0.49}}
\multiput(64.15,15.12)(0.06,0.5){1}{\line(0,1){0.5}}
\multiput(64.22,15.62)(0.05,0.5){1}{\line(0,1){0.5}}
\multiput(64.27,16.11)(0.04,0.5){1}{\line(0,1){0.5}}
\multiput(64.31,16.61)(0.03,0.5){1}{\line(0,1){0.5}}
\multiput(64.34,17.11)(0.02,0.5){1}{\line(0,1){0.5}}
\multiput(64.36,17.61)(0,0.5){1}{\line(0,1){0.5}}
\multiput(64.35,18.61)(0.01,-0.5){1}{\line(0,-1){0.5}}
\multiput(64.33,19.11)(0.02,-0.5){1}{\line(0,-1){0.5}}
\multiput(64.3,19.61)(0.03,-0.5){1}{\line(0,-1){0.5}}
\multiput(64.26,20.11)(0.04,-0.5){1}{\line(0,-1){0.5}}
\multiput(64.2,20.61)(0.06,-0.5){1}{\line(0,-1){0.5}}
\multiput(64.13,21.1)(0.07,-0.5){1}{\line(0,-1){0.5}}
\multiput(64.05,21.6)(0.08,-0.49){1}{\line(0,-1){0.49}}
\multiput(63.96,22.09)(0.09,-0.49){1}{\line(0,-1){0.49}}
\multiput(63.86,22.58)(0.1,-0.49){1}{\line(0,-1){0.49}}
\multiput(63.74,23.06)(0.12,-0.49){1}{\line(0,-1){0.49}}
\multiput(63.62,23.55)(0.13,-0.48){1}{\line(0,-1){0.48}}
\multiput(63.48,24.03)(0.14,-0.48){1}{\line(0,-1){0.48}}
\multiput(63.33,24.51)(0.15,-0.48){1}{\line(0,-1){0.48}}
\multiput(63.16,24.98)(0.16,-0.47){1}{\line(0,-1){0.47}}
\multiput(62.99,25.45)(0.17,-0.47){1}{\line(0,-1){0.47}}
\multiput(62.81,25.91)(0.09,-0.23){2}{\line(0,-1){0.23}}
\multiput(62.61,26.37)(0.1,-0.23){2}{\line(0,-1){0.23}}
\multiput(62.4,26.83)(0.1,-0.23){2}{\line(0,-1){0.23}}
\multiput(62.18,27.28)(0.11,-0.23){2}{\line(0,-1){0.23}}
\multiput(61.96,27.72)(0.11,-0.22){2}{\line(0,-1){0.22}}
\multiput(61.72,28.16)(0.12,-0.22){2}{\line(0,-1){0.22}}
\multiput(61.47,28.6)(0.13,-0.22){2}{\line(0,-1){0.22}}
\multiput(61.21,29.02)(0.13,-0.21){2}{\line(0,-1){0.21}}
\multiput(60.93,29.44)(0.14,-0.21){2}{\line(0,-1){0.21}}
\multiput(60.65,29.86)(0.14,-0.21){2}{\line(0,-1){0.21}}
\multiput(60.36,30.26)(0.15,-0.2){2}{\line(0,-1){0.2}}
\multiput(60.06,30.66)(0.1,-0.13){3}{\line(0,-1){0.13}}
\multiput(59.75,31.06)(0.1,-0.13){3}{\line(0,-1){0.13}}
\multiput(59.43,31.44)(0.11,-0.13){3}{\line(0,-1){0.13}}
\multiput(59.1,31.82)(0.11,-0.13){3}{\line(0,-1){0.13}}
\multiput(58.77,32.19)(0.11,-0.12){3}{\line(0,-1){0.12}}
\multiput(58.42,32.55)(0.12,-0.12){3}{\line(0,-1){0.12}}
\multiput(58.06,32.9)(0.12,-0.12){3}{\line(1,0){0.12}}
\multiput(57.7,33.24)(0.12,-0.11){3}{\line(1,0){0.12}}
\multiput(57.33,33.58)(0.12,-0.11){3}{\line(1,0){0.12}}
\multiput(56.95,33.9)(0.13,-0.11){3}{\line(1,0){0.13}}
\multiput(56.56,34.22)(0.13,-0.11){3}{\line(1,0){0.13}}

\linethickness{0.1mm}
\multiput(64.53,39.84)(0.18,0.12){125}{\line(1,0){0.18}}
\linethickness{0.1mm}
\put(72.03,9.84){\line(1,0){30}}
\linethickness{0.1mm}
\multiput(71.56,44.06)(0.12,-0.49){70}{\line(0,-1){0.49}}
\linethickness{0.1mm}
\multiput(19.53,9.84)(0.48,0.12){63}{\line(1,0){0.48}}
\linethickness{0.1mm}
\multiput(87.03,32.34)(0.48,0.12){63}{\line(1,0){0.48}}
\linethickness{0.5mm}
\put(74.38,32.34){\line(1,0){12.65}}
\linethickness{0.5mm}
\multiput(74.38,32.34)(0.12,-0.44){16}{\line(0,-1){0.44}}
\linethickness{0.5mm}
\multiput(76.25,25.31)(0.18,0.12){59}{\line(1,0){0.18}}
\put(18.59,7.5){\makebox(0,0)[cc]{$\scriptstyle 0$}}

\put(49.06,14.53){\makebox(0,0)[cc]{$\scriptstyle x$}}

\put(72.03,7.03){\makebox(0,0)[cc]{$\scriptstyle v_1$}}

\put(59.84,41.25){\makebox(0,0)[cc]{$\scriptstyle v_2$}}

\put(83.28,7.03){\makebox(0,0)[cc]{$\scriptstyle(1+\varepsilon/2)v_1$}}

\put(64.06,45.47){\makebox(0,0)[cc]{$\scriptstyle (1+\varepsilon/2)v_2$}}

\put(60.78,7.03){\makebox(0,0)[cc]{$\scriptstyle (1-\varepsilon)v_1$}}

\put(49.53,35.16){\makebox(0,0)[cc]{$\scriptstyle (1-\varepsilon)v_2$}}

\put(95.47,30.47){\makebox(0,0)[cc]{$\scriptstyle v_1+v_2-x$}}

\put(76.25,16.88){\makebox(0,0)[cc]{$\scriptstyle C_n$}}

\put(84.16,24.31){\makebox(0,0)[cc]{$\scriptstyle \Delta^\star$}}

\put(71.09,45.94){\makebox(0,0)[cc]{}}

\put(67.81,36.09){\makebox(0,0)[cc]{}}

\put(79.06,24.84){\makebox(0,0)[cc]{}}

\put(0,35){\makebox(0,0)[cc]{}}

\put(123,41){\makebox(0,0)[cc]{$\scriptstyle v_1+v_2$}}

\end{picture}

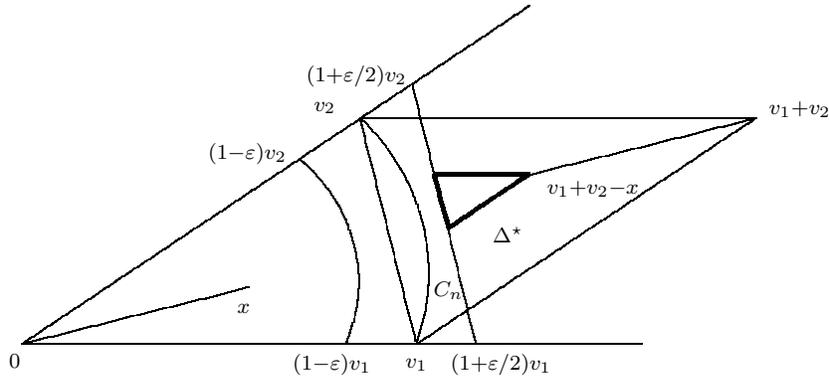
\centerline{Figure 2. The proof of Lemma 7.2}

\section{Proof of Theorem~\ref{th:onC}}\label{sec:proofs}

In this section we prove Theorem~\ref{th:onC} apart from the uniqueness of $C$ and $r$ which
will be shown in the next section.

The Blaschke selection theorem and Corollary~\ref{cor:sandw} imply that every subsequence of
$(\area C_n)^{-1/2}C_n$ contains a convergent (in Hausdorff metric) subsequence.
Corollary~\ref{cor:liminf} guarantees then the existence of positive integers $n_1<n_2<\dots$
such that $\lim n_k^{-3/2}L_{n_k}=\al$ and $\lim \dist((\area C_{n_k})^{-1/2}C_{n_k},C)=0$
for some convex body $C \in \ck$. Define $\la_k=\sqrt{\area C_{n_k}}$ and set, for simpler
writing, $C^k=\la_k^{-1}C_{n_k}$. It is evident that $rB \subset C \subset RB$, showing that,
for every $\de>0$, $(1-\de)C \subset C^k \subset (1+\de)C$ for all large enough $k$. Since
$n_k=\frac 6{\pi^2}\area C_{n_k}(1+o(1))$, $\la_k=\frac {\pi}{\sqrt 6}\sqrt {n_k}(1+o(1))$.

It follows immediately that $\area C=1$. We show next that $\int_C zdz =0$. For this it
suffices to prove that $\int_Cf(z)dz=0$ in the case when $f$ is the linear function $f(z)=x$
and $f(z)=y$ where $z=(x,y)$. Choose $\eps>0$ and then, using Lemma~\ref{le:approx}, $k_0$ so
large that, for $k>k_0$,
\[
(1-\eps/2)C_{n_k} \cap \p \subset E_{n_k} \subset C_{n_k} \cap \p.
\]
It follows now that there is a $k_1$ so that for all $k>k_1$
\begin{equation}\label{eq:sandw}
(1-\eps)\la_kC \cap \p \subset E_{n_k} \subset (1+\eps)\la_kC \cap \p.
\end{equation}
Using the notation $f(X)=\sum_{z \in X} f(z)$ when $X\subset \r$ is finite, we have
$f(E_{n_k})=0$. Next,
\begin{eqnarray*}
|f(\p\cap \la_kC)|&=&|f(\p\cap \la_kC)-f(E_{n_k})|\\ &\le& |f\left(\p\cap
[(1+\eps)\la_kC\setminus (1-\eps)\la_kC]\right)|\\ &\ll& \eps \la_k \max \{f(z): z\in \la_kC\}
\ll \eps n_k.
\end{eqnarray*}
On the other hand, by (\ref{eq:fpla}),
\[
|f(\p\cap \la_kC)|=\frac 6{\pi^2} \la_k^3 \int_Cf(z)dz\left(1+O(\la_k^{-1}\log \la_k)\right)
\]
as one can check easily. So if $\int_Cf(z)dz\ne 0$, then $f(\p\cap \la_kC)$ is of order
$n_k^{3/2}$. But as we have just shown, $|f(\p\cap \la_kC)| \ll \eps n_k$. So indeed,
$\int_Cf(z)dz= 0$, or, in other words, $g(C)=0$.

An almost identical proof, this time with the 1-homogeneous function $f(z)=\|z\|$ gives
$$\frac {\pi}{\sqrt 6} \int_C\|x\|dx=\al(D).$$
We only give a sketch: Equation (\ref{eq:sandw}) shows that
$$\left|\sum_{z \in \p\cap \la_kC}\|z\|-\sum_{z \in \p\cap E_{n_k}}\|z\|\right| \ll \eps n_k.$$
Here $\sum_{z \in \p\cap E_{n_k}}\|z\|=L_{n_k}$ and so $\lim n_k^{-3/2}\sum_{z \in \p\cap \la_kC}\|z\|=\al(D)$.
The estimate (\ref{eq:fp}) says now that
\[
\left|\sum_{z \in \p\cap \la_kC}\|z\|- \frac 6{\pi^2}\int_{\la_kC}\|x\|dx\right| \ll n_k \log n_k,
\]
and $\frac {\pi}{\sqrt 6} \int_C\|x\|dx=\al(D)$ follows.

Lemma~\ref{l:conn} applies now because $g(C)=0$ and $\area C=1$. So there is a sequence $Q_n
\in \pn$ with $\lim n^{-3/2}\per Q_n =\al(D)$. Then $L_n \le \per Q_n$ implies that $\lim
n^{-3/2}L_n =\al(D)$. \hfill $\Box$\medskip

\section{The variational problem}

Next we turn to uniqueness. As first step we treat a special case.

\begin{lemma}\label{centered}
Let $r_0$ be the radial function of $D \in \ck$ with $g(D)=0$. Then $r_0$ is the unique
solution to $VP(r_0)$.
\end{lemma}

{\bf Proof.} We consider the variational problem which ignores the constraints about the
center of gravity :

$$\mbox{ minimize }\int_0^{2\pi}  r^3(t)/r_0 (t)dt $$
$$\mbox{subject to}\quad\int_0^{2\pi} r^2(t)dt =2$$

From H\"older's inequality :
$$  \int_0^{2\pi}r^2\leq \left(\int_0^{2\pi}{r^3\over r_0}\right)^{2/3}\left(\int_0^{2\pi}r_0^2\right)^{1/3}           $$
which is an equality if and only if $r$ and $r_0$ are proportional. In our case
$\int_0^{2\pi}r^2=\int_0^{2\pi}r_0^2=2$ and so $r=r_0$. \hfill $\Box$

\medskip

We now use the previous lemma to treat the general case:

\begin{lemma}
There exists a unique solution $r\in \cf$ to problem $VP(r_0)$. This solution is equal to
$$r=\left({a\over r_0}+b\cos t +c\sin t\right)^{-1}    $$
where $a>0,\;b,c$ are the unique real numbers which make the function $r$ satisfy the three
constraints of $VP(r_0)$.
\end{lemma}

{\bf Proof.} We prove in Appendix 3 that every optimal solution $r \in \cf^c$ to $VP(r_0)$ is
of the form $r(t)=({a\over r_0}+b\cos t+c\sin t)^{-1}$ with suitable constants $a,b,c \in R$.
We have shown that the radial function, $r(t)$, of $C$ from Theorem~\ref{th:onC} is an
optimal solution to $VP(r_0)$. As $C$ is convex, $r(t)$ is equal to $({a\over r_0}+b\cos
t+c\sin t)^{-1}$. According to the previous Lemma, the unique solution to the variational
problem $VP(r)$ is $r$.

Consider now another optimal solution, $r^*$, to $VP(r_0)$. It is clear that $r^*$ is a
feasible solution to $VP(r)$ and that
$$\int_0^{2\pi}\frac {r^{*3}}{r_0}=\int_0^{2\pi}\frac {r^{3}}{r_0}    .$$
Further,
\[
a\int_0^{2\pi}\frac {r^{*3}}{r_0}=\int_0^{2\pi}r^{*3}\left(\frac a{r_0}+b\cos t +c\sin t\right)=\int_0^{2\pi}\frac {r^{*3}}{r},
\]
and, in the same way,
\[
a\int_0^{2\pi}\frac {r^{3}}{r_0}=\int_0^{2\pi}r^{3}\left(\frac a{r_0}+b\cos t +c\sin t\right)=\int_0^{2\pi}\frac {r^{3}}{r}.
\]

So $r^*$, too, is an optimal solution to $VP(r)$. By the Lemma, $r=r^*$, and $a>0$ follows as
well. \hfill$\Box$

{\bf Remark:} After reading this proof, one easily understands  that $r(t)$ is the 
unique solution to the variational problem in a class of functions larger than $\cf$.

\section{Proof of Theorem~\ref{th:onPn}}\label{sec:PndistP}

This is fairly simple once we know that $C$ is unique. Let $u(t)=(\cos t,\sin t)$ be the unit
vector in direction $t \in [0,2\pi]$. When a minimizer $P_n$ is translated as
Theorem~\ref{th:onPn} specifies, the sum of the edges of $P_n$ having direction between
$u(0)$ and $u(t)$ is very close to the sum of the primitive vectors having direction between
$u(0)$ and $u(t)$ in $C_n$. The latter, divided by $n^{3/2}$ is very close to
$P(t)=\int_{C(t)}zdz$ where $C(t)$ is the set of vectors in $C$ with direction between $u(0)$
and $u(t)$. The curve $P(t)$ is closed (because $g(C)=0$) and convex (this has been shown in
\cite{bapr}), so it is the boundary of a convex set $P$. The simple and straightforward
checking of
\[
\lim \dist (n^{-3/2}P_n,P)=0
\]
is left to the reader. We remark that the convexity of $P(t)$ follows also from the fact that
the boundary of $P_n$, after suitable rescaling, tends to $P(t)$. \hfill$\Box$
\medskip

The same construction $C \to P$ with $P(t)=\int_{C(t)}zdz$ is used, with a similar purpose,
in \cite{bapr}. Further properties of the construction are also established there.

\section{An example}

We concentrate now on the cases when the solution is constant which correspond to the case
when the limit shape of the polygon is a circle.

\begin{lemma}
The solution is constant if and only if $1/r_0$ is of the form $a+b\cos\theta+c\sin\theta$,
or, in other words, when $r_0$ is the radial function of an ellipse having its focus point at
the origin.
 \end{lemma}

{\bf Proof.} Suppose the solution is constant, the form of $r_0$ is then directly derived
from Lemma 9.2. Conversely, if ${1\over r_0}$ is of the form $a+b\cos\theta+c\sin\theta$, the
solution is then also of the form ${1\over r}=a'+b'\cos\theta+c'\sin\theta$. This says that
it is the radial function of an ellipse having its focus point at the origin. We conclude by
observing that the only ellipses whose centre of gravity is at the same time their focus
point, are circles.  \hfill$\Box$

\section{Appendix 1}

\begin{lemma}
Let $M_n(w)$ be the sum of the lengths of the $n$ shortest (in Euclidean norm) distinct
vectors in $\z$ lying in a strip of width $w$, centred at the origin. Suppose $\ga \in
(0,1/2]$, then $w \le \ga n^{1/2}$ implies $M_n(w) \gg n^{3/2}/\ga $.

\end{lemma}

{\bf Proof.}  It is clear that this set of vectors is just the set of lattice points contained
in $A:=dB \cap T$ where $T$ is a strip of width $w$, centred at the origin, and $d$ is a
suitable radius making $A\cap \z$ have exactly $n$ elements (ties broken arbitrarily). Let
$\phi$ denote the angle that the strip $T$ makes with the $x$-axis of $\r$. We may assume by
symmetry that $\phi \in [0,\pi/4]$.

Observe first that $d \ge \sqrt n /2$ since otherwise the disk $dB$ would contain fewer than
$n$ lattice points. Let $Q(z)$ denote the unit square centred at $z \in \r$ and let $\ell_k$
be the line with equation $x=k$ ($k$ is an integer). Clearly, $\ell_k$ intersects $S$ in a
segment of length $w\cos \phi$, and so $\ell_k \cap \z$ contains at least $\lfloor w/\cos
\phi \rfloor$ and at most $\lfloor w/\cos \phi \rfloor +1$ lattice points from $S$.

Assume first that $w/\cos \phi \ge 1$. As is easy to see, $\area A \cap Q(z)$ is at least
$1/4$ for $z \in A \cap \z$.   Hence, $\area A \ge n/4$. Since $\area A < 2dw$, $d > n/(4w)$
follows.

For simpler notation write $u=(d\cos\phi)/2$. For the lines $\ell_k$ with $k \in [u,2u-w/2]$,
$\ell_k \cap A$ contains at least $\lfloor w/\cos \phi \rfloor$ lattice points. Since $w<u$, 
there are at least $\lfloor 2u -w/2 \rfloor-\lfloor u \rfloor\gg u$ such lines. All of them
have distance at least $(d-w)/2\gg d$ from the origin. Consequently, using the bounds $w\le
\ga n^{1/2}$ and $d\ge n^{1/2}/2$ generously,
\begin{eqnarray*}
M_n(w) \gg d \left\lfloor \frac{w}{\cos \phi}\right\rfloor u  \gg
d^2w \ge \left(\frac {n}{4w}\right)^2 w \gg \frac 1{\ga} n^{3/2}.
\end{eqnarray*}

Assume next that $w/\cos \phi <1$. There are at most six $z \in A \cap \z$ such that $Q(z)$ intersects the
boundary of $dB$. For the other $z \in A \cap \z$, $Q(z)$ intersects the boundary of $A$ in
one or two line segments, whose total length is between $1/\cos \phi$ and $2/\cos \phi$. For
distinct lattice points in $A \cap \z$ the corresponding segments do not overlap. This
implies that
\[
\frac {n-6}{\cos \phi} \le 4d \le \frac {2(n-6)}{\cos \phi}.
\]
Each line $\ell_k$ with $|k| \le 2u/3$ contains at most one lattice point from $A$. The
remaining points from $A\cap \z$, and there are at least $n-2\lfloor 2u/3 \rfloor-1$ of
them, are at distance $\frac d3-1$ from the origin. Hence, we see
\[
M_n(w) \ge \left(\frac d3 -1\right)\left(n-\lfloor 2\frac{d\cos \phi}3 \rfloor-1\right) \gg
n^2.
\]
\hfill$\Box$

\section{Appendix 2}

We start the proof of Lemma~\ref{le:approx} with the following Claim.

\begin{cl} Suppose $a,b,c,d$ are vertices of $K$ (in anticlockwise order), $[a,b]$ and $[c,d]$ are
edges of $K$, and $\angle b0c < 3\de$. Let $x$ be the intersection point of the lines through
$a,b$ and $c,d$, and let $y$ be the intersection point of the lines through $0,x$ and $a,c$.
Then $|x-y|\le 4\de (R/r)^2 |y|$.
\end{cl}

{\bf Proof.} The condition $rB \subset K \subset RB$ implies that $\be =\angle 0xb = \angle
0ba - \angle xba > \arcsin r/R - 3\de$ since
$$\sin \angle 0ba={d(0, \ell_{a,b})\over |b|}$$
($\ell_{a,b}$ being the line through $a$ and $b$) $|b|<R$ , $d(0, \ell_{a,b})>r$ by
assumption, so that $\sin \angle 0ba>r/R$, see Figure 3.

 Further $\angle xyc =\angle 0xa-\angle x0b
>\be$. The sine theorem in the triangle $x,y,c$ shows that
\[
\frac {|x-y|}{|x-c|}=\frac {\sin \angle cxy}{\sin \angle cyx},
\]
and similarly, the sine theorem in the triangle $x,0,c$ shows that
\[
\frac {|x-c|}{|x|}=\frac {\sin \angle c0x}{\sin \angle 0cx}.
\]
Multiplying them gives
\[
\frac {|x-y|}{|x|}=\frac {\sin \angle cxy \sin \angle c0x}{\sin \angle cyx \sin \angle 0cx}<
\frac {\sin 3\de}{(r/R)\sin \be}.
\]
Next, since $|y|=|x|-|x-y|$, we have
\[
\frac {|x|}{|x|-|x-y|}=\frac {1}{1- \frac {|x-y|}{|x|}}< \frac {1}{1-\frac {\sin
3\de}{(r/R)\sin \be}}
\]
We use this inequality next in the form
\[
\frac {|x-y|}{|y|}<\frac {\sin 3\de}{(r/R)\sin \be}\cdot \frac {|x|}{|x|-|x-y|}<\frac {\sin
3\de}{(r/R)\sin \be -\sin 3\de}< 4\de \left(\frac Rr\right)^2,
\]
where we only have to check the validity of the last inequality. This is a matter of direct
computation using that $\sin \be > \sin (\arcsin (r/R)-3\de) > (r/R)\cos 3\de- \sin 3\de$ and
the assumption that $\de < 0.02(r/R)^2$ implying, in particular, that $\de < 0.02$. What is
to be checked now is that
\[
\tan 3\de\left[1+4\de\left(\frac Rr\right)^2\left(\frac rR +1\right)\right]\le 4\de.
\]
Here $\de (R/r)^2 < 0.02$ and so the expression in the square bracket is at most $1.16$ and the inequality follows.
We omit
the details. \hfill $\Box$
\bigskip

\ifx\JPicScale\undefined\def\JPicScale{1}\fi
\unitlength \JPicScale mm
\begin{picture}(102.88,74.94)(0,0)
\linethickness{0.3mm}
\qbezier(32.47,29.5)(32.55,30.42)(35.53,39.1)
\qbezier(35.53,39.1)(38.51,47.78)(43.3,53.15)
\qbezier(43.3,53.15)(50.06,59.55)(58.52,62.85)
\qbezier(58.52,62.85)(66.99,66.14)(75.8,64.98)
\qbezier(75.8,64.98)(82.48,63.91)(88.12,59.67)
\qbezier(88.12,59.67)(93.76,55.43)(97.47,49.21)
\qbezier(97.47,49.21)(100.47,43.98)(101.22,37.86)
\qbezier(101.22,37.86)(101.97,31.74)(101.08,25.56)
\qbezier(101.08,25.56)(100.08,20.65)(97.16,15.43)
\qbezier(97.16,15.43)(94.23,10.22)(93.85,9.79)
\linethickness{0.3mm}
\multiput(57.74,9.79)(0.12,0.46){120}{\line(0,1){0.46}}
\linethickness{0.3mm}
\multiput(57.74,9.79)(0.12,0.12){328}{\line(1,0){0.12}}
\linethickness{0.3mm}
\multiput(57.74,9.79)(0.12,0.21){301}{\line(0,1){0.21}}
\linethickness{0.3mm}
\multiput(93.85,72.86)(0.12,-0.86){71}{\line(0,-1){0.86}}
\linethickness{0.3mm}
\multiput(32.92,51.18)(0.34,0.12){181}{\line(1,0){0.34}}
\put(55.55,8.35){\makebox(0,0)[cc]{$\scriptstyle 0$}}

\put(102.88,22.11){\makebox(0,0)[cc]{$\scriptstyle a$}}

\put(100.62,50.2){\makebox(0,0)[cc]{$\scriptstyle b$}}

\put(97.01,74.34){\makebox(0,0)[cc]{$\scriptstyle x$}}

\put(92.05,65.96){\makebox(0,0)[cc]{$\scriptstyle\beta$}}

\put(70.83,69.41){\makebox(0,0)[cc]{$\scriptstyle c$}}

\put(43.75,58.08){\makebox(0,0)[cc]{$\scriptstyle d$}}

\linethickness{0.3mm}
\multiput(72.19,64.98)(0.12,-0.17){241}{\line(0,-1){0.17}}
\put(85.28,51.67){\makebox(0,0)[cc]{$\scriptstyle y$}}

\put(34.62,74.94){\makebox(0,0)[cc]{}}

\end{picture}

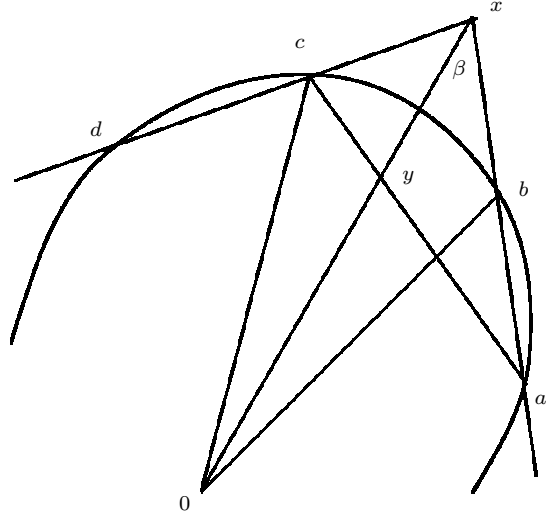
\centerline{Figure 3. The proof of Claim 13.1} 

\medskip

The {\bf Proof} of Lemma~\ref{le:approx} is an algorithm that constructs the vertex set $V$
of $Q$. We start with $V=\emptyset$. We call the edge $[a,b]$ of $K$ {\sl special} if $\angle
a0b \ge \de$. Let $W$ be a cone with apex at $0$ and angle $\de$. It follows that if $W$ is
disjoint from all special edges, then it contains a vertex of $K$.

{\bf Case 1.} Let $[a_1,b_1],[a_2,b_2],\dots,[a_k,b_k]$ be consecutive special edges in
anticlockwise order so that $\angle b_i0a_{i+1} <3\de$ for all $i=1,\dots,k-1$ (or up to $k$
if $\angle b_k0a_1 < 3\de$). We call this a maximal chain of consecutive special edges if
there is no special edge $[a,b]$ with $\angle b0a_1 < 3\de$ or $\angle b_k0a < 3\de$.

For such a maximal chain we put the vertices $a_1,\dots,a_k,b_k$ (or $a_1,\dots,a_k$ if
$\angle b_k0a_1 < 3\de$) into $V$, and we do so for all such maximal chains.

{\bf Case 2.} Let $[a_1,b_1]$ and $[a_2,b_2]$ be consecutive special edges with vertices
$a_1,b_1,a_2,b_2$ in anticlockwise order so that $\ga:= \angle b_10a_2 \ge 3\de$. Then we
choose $\de'\in [\de, 3\de]$ so that $\ga/\de'$ is an odd integer, say $2h+1$. This is always
possible since there is an odd integer between $\ga/(3\de)$ and $\ga/\de$ because their
difference is $\ga/\de -\ga/(3\de)=2\ga/(3\de)\ge 2$.

Subdivide now the cone $\pos\{b_1,a_2\}$ into $2h+1$ subcones, each of angle $\de'$ and pick
a vertex $u_1,\dots,u_h$ from every second subcone. Finally, put $b_1,u_1,\dots,u_h,a_2$ into
$V$.

If there are only two special edges $[a_1,b_1]$ and $[a_2,b_2]$, then one has to do the same
construction between edges $[a_2,b_2]$ and $[a_1,b_1]$ as well. If there is only one special
edge, then the construction is carried out from $b_1$ to $a_1$ as if one had two special
edges $[a_1,b_1]$ and $[b_1,a_1]$.

Finally, if there are no special edges, then we chose a $\de' \in [\de,2\de]$ so that
$2\pi/\de'$ is an even integer, $2h$, say. This is evidently possible. Subdivide the plane
into cones of angle $\de'$ (with apex at $0$) and choose a vertex $u_1,\dots,u_h$ from every
second cone, and set $V=\{u_1,\dots,u_h\}$.

The algorithm is finished. By construction $\angle v_i0v_{i+1} \ge \de$: for the angle at
$0$. Finally we check condition $K \subset (1+4\de(R/r)^2)Q$. Let
$v_i,v_{i+1},v_{i+2},v_{i+3}$ be four consecutive vertices of $Q$ in anticlockwise order.
Rename these points as $a,b,c,d$ as in the Claim. Then $K \cap \pos (b,c) \setminus Q$ is
contained in the triangle $b,c,x$ from the Claim. Now $y \in Q$ because $y$ lies on the
segment $[a,c]$, and so $x \in (1+4\de(R/r)^2)Q$ according to the Claim. So the triangle
$b,c,x$ is contained $(1+4\de(R/r)^2)Q$. \hfill $\Box$
\medskip

\section{Appendix 3}

It happens that standard theorems of the Calculus of Variations (see for instance \cite{Sagan}) are stated in a  $C^1$
setting, and suppose also that  the function  $r_0$ involved in the problem is $C^1$. Since
these conditions are not satisfied in our problem, we have to elaborate the following
statement:

\begin{lemma}
All the  solutions $r\in \cf^c$ satisfying problem $VP(r_0)$ are of the form
$$r=\left({a\over r_0}+b\cos t +c\sin t\right)^{-1}    $$
where $a,b,c$ are real numbers which make the function $r$ satisfy the three
constraints of problem (\ref{eq:var}).
\end{lemma}

{\bf Proof.} Consider $r$ an optimal solution in $\cf^c$. Let $h$ be a function on $[0,
2\pi]$ such that the   perturbed function  $r_\eps:=r+\eps h$ remains in $\cf^c$ for $\eps$
in a neighbourhood of $0$ (notice that all the twice differentiable functions are
convenient). This perturbation won't be feasible in general. We want  to modify it in order
to make it feasible. That is what we do in the two first steps.

{\bf Step 1.} We translate the set defined by the function $r_\eps$ in order  to get a
centred set defined by a new radial function $\tilde r_\eps$ we evaluate up to some
$o(\eps)$.

In the following, the notation $o(\eps)$ stands for some family of functions, which may be constant, indexed by $\eps$, such
that both ${o(\eps)\over\eps}$ converges to $0$ as $\eps$ goes to $0$, and
${o(\eps)\over\eps}$ is dominated.

The coordinates of the centre of gravity of the set defined by $r_\eps$ are
$$\left(\int_0^{2\pi}  r_\eps^3\cos t dt,    \int_0^{2\pi}  r_\eps^3 \sin t dt\right)= \eps \left(\int_0^{2\pi}  3r^2h\cos t dt,\int_0^{2\pi}  3r^2h\sin t dt\right)+o(\eps)$$

Recall that $u(t)=(\cos t,\sin t)$. Define the numbers $r_h,\theta_h$ by setting $
r_hu(\theta_h) :=(\int_0^{2\pi} 3r^2h\cos t dt,\int_0^{2\pi} 3r^2h\sin t dt).$

For a given $\theta$, the polar coordinates of $r_\eps u(\theta)-\eps r_hu(\theta_h)+o(\eps)$
are given by
$$\tilde\theta(\theta):= \theta-\eps r_h {\sin (\theta_h-\theta)\over r(\theta)}+o(\eps)$$
$$\tilde r(\theta)=r_\eps(\theta)-\eps r_h\cos (\theta_h-\theta)+o(\eps)$$
Hence, $\tilde r$ can be expressed as a function of $\tilde\theta$ as follows:
$$\tilde r(\tilde\theta)=r_\eps\left(\tilde\theta+\eps r_h {\sin (\theta_h-\tilde\theta)\over r(\theta)}+o(\eps)\right)-\eps r_h\cos \left(\theta_h-\tilde\theta+o(\eps)\right)+o(\eps)$$
Using now the almost everywhere differentiability of $r$ (and therefore of $r_\eps$) which is
inherited from convexity, we obtain that, almost everywhere,
$$\tilde r(\theta)=r(\theta)+\eps\left[ h(\theta)+ r_h \sin (\theta_h-\theta){r'\over r}(\theta)- r_h\cos(\theta_h-\theta)\right]+o(\eps)$$
Note that the domination of ${o(\eps)\over\eps}$ in the last step is due to the fact that the
left and right derivatives of $r$ are bounded on $[0,2\pi]$.

{\bf Step 2.} We obtain a completely feasible function $r^f$, by normalizing $\tilde r$ by
the area of the set defined by $\tilde r$, which is the same as the area of the set defined
by $r_\eps$, since the two sets  are obtained one from the other  by a translation.

Define, $$  r^f(\theta)={\tilde r(\theta)\over({{1\over2}\int_0^{2\pi}r_\eps^2)^{1/2}}}=
  {\tilde r(\theta)\over({{1\over2}\int_0^{2\pi}(r+\eps h)^2)^{1/2}}}=
  \tilde r (\theta)\left(1-{\eps\over2}(\int_0^{2\pi}rh)+o(\eps)\right)$$
The function  $r^f(\theta)$ can be written as $r(\theta)$ times the function
$$1+\eps\left[ {h(\theta)\over r(\theta)}+r_h \sin (\theta_h-\theta){r'(\theta)\over r^2(\theta)}- r_h{\cos(\theta_h-\theta)\over r(\theta)}-{1\over2}(\int_0^{2\pi}rh)\right]+o(\eps)$$

{\bf Step 3.} Now, we test the optimality of the function $r$ by considering the functional
applied to the feasible perturbation $r^f$ and writing the integral $\int_0^{2\pi}
{(r^f)^3\over r_0}$ as $\int_0^{2\pi} {r^3\over r_0}$ plus

$$3\eps\int_0^{2\pi} {r^3\over r_0}\left[{h(\theta)\over r(\theta)}(\theta)+ r_h \sin (\theta_h-\theta){r'(\theta)\over r^2(\theta)}(\theta)- r_h{\cos(\theta_h-\theta)\over r(\theta)}-{1\over2}(\int_0^{2\pi}rh)\right]+o(\eps)$$

When developing the sine and cosine in the above bracket and performing the integration on
$\theta$ (and keeping in mind that $r_h$ and $\theta_h$ are constants that don't depend on
$\theta$ !) we deduce that, if $r$ is optimal,  there exist real constants $A$, $B$ and $C$
such that for all twice differentiable function $h$,
 $$\int_0^{2\pi}{r^2 h\over r_0}+A r_h\sin\theta_h+B r_h\cos\theta_h+C\int_0^{2\pi}rh=0$$
Recall that $(r_h\cos\theta_h, r_h\sin\theta_h)=      (\int_0^{2\pi}  3r^2h\cos t dt,\int_0^{2\pi}  3r^2h\sin t dt)   $.

Therefore, for all twice differentiable functions $h$
$$   \int_0^{2\pi}hr^2\left({1\over r_0}+3A \cos\theta+3B \sin\theta+{C\over r}\right)=0   $$
which implies that the bracket inside the integral is 0.
\hfill$\Box$

\section{Acknowledgements} The first author was partially supported by Hungarian National
Science Foundation Grants T 032452 and T 60427, and also by the Discrete and Convex Geometry
project, MTKD-CT-2005-014333, of the European Community. The second author was partially
supported by ANR grant MEMEMO.

Imre B\'ar\'any \\
R\'enyi Institute of Mathematics,\\
Hungarian Academy of Sciences\\
H-1364 Budapest Pf.~127  Hungary\\
{\tt barany@renyi.hu}\\
and\\
Department of Mathematics\\
University College London\\
Gower Street, London, WC1E 6BT, UK
\and\bigskip
\\
Nathana\"el Enriquez \\
Laboratoire Modal'X\\
Universit\'e Paris-Ouest\\
200 Avenue de la R\'epublique, 92001 Nanterre, France\\
{\tt nenriquez@u-paris10.fr} \\and\\
Laboratoire de Probabilit\'e et Mod\`eles Al\'eatoires \\
Universit\'e Paris 6 \\
4 Place Jussieu, 75005 Paris, France 

\end{document}